
\magnification = \magstep1
\baselineskip = 1.5\baselineskip
\raggedbottom

\font\title = cmbx12

\def\sqr#1#2{{\vcenter{\vbox{\hrule height.#2pt
\hbox{\vrule width.#2pt height#1pt \kern#1pt
\vrule width.#2pt}
\hrule height.#2pt}}}}

\def\square
{\mathchoice\sqr44\sqr44\sqr{6}8\sqr{4}8}

\def\qed{$\quad \square$}

\def\newpar{\vskip 0.25\baselineskip}

\def\state#1 {{\narrower\narrower\smallskip \noindent {\rm #1} \smallskip
\par}}

\def\result#1#2 {{\narrower\smallskip \noindent {\bf #1:} {\sl #2} \smallskip
\par}}

\def\ref#1#2 {\item{[#1]} {#2}}


\centerline {\title On the Negative Case of the}
\newpar
\centerline {\title Singular Yamabe Problem}

\vskip \baselineskip
\centerline {David L.\ Finn}
\centerline {Merrimack College}

\vfill
\noindent
{\bf Abstract:} Let $(M,g)$ be a compact Riemannian manifold of dimension $n
\geq 3$, and let $\Gamma$ be a nonempty closed subset of $M$.  The negative
case of the Singular Yamabe Problem concerns the existence and behavior of a
complete metric $\hat g$ on $M\backslash\Gamma$ that has constant negative
scalar curvature and is pointwise conformally related to the smooth metric $g$.
Previous results have shown that when $\Gamma$ is a smooth submanifold
of dimension $d$ there exists such a metric if and only if $d > {n-2 \over 2}$.
In this paper, we consider a general class of closed sets and show the
existence of a complete conformal metric  $\hat g$ with constant negative
scalar curvature depends on the dimension of the {\it tangent cone to $\Gamma$}
at every point.  Specifically, provided $\Gamma$ admits a nice tangent cone at
$p$, we show that when the dimension of the tangent cone to $\Gamma$ at $p$ is
less than ${n-2 \over 2}$ then there can not exist a {\it negative Singular
Yamabe metric} $\hat g$ on $M\backslash\Gamma$.

\vfill
{{\narrower {

\newpar
}\par}}

\vfill\eject

\line {\bf 1. Introduction: \hfil}

\newpar
The resolution of the {\it Yamabe Problem} by Schoen (cf.\ [LP] or [Sc1])
showed that every compact Riemannian manifold of dimension $n \geq 3$ is
conformally equivalent to one with constant scalar curvature.  A natural
question to then ask is whether every noncompact Riemannian manifold of
dimension $n \geq 3$ is conformally equivalent to a complete manifold with
constant scalar curvature.  For noncompact manifolds with a simple structure at
infinity, this question may be studied by solving the so-called {\it Singular
Yamabe problem}:

\state {Given a compact Riemannian manifold $(M,g)$ of dimension $n \geq 3$,
and a nonempty closed set $\Gamma$ in $M$.  Find a complete metric $\hat g$ on
$M\backslash\Gamma$ with constant scalar curvature that is pointwise
conformally related to $g$.}

\noindent
The principal benefit of studying this version of the Yamabe Problem on
noncompact manifolds is that the complete metric $\hat g$ has a uniform
underlying structure at infinity.  The general problem of conformally
deforming an arbitrary metric $\bar g$ to a complete metric $\hat g$ with
constant scalar curvature requires extra conditions on the metric structure of
$\bar g$ at infinity (for example see [LTY]).

\newpar
In this paper, we consider the {\it negative case}, when the complete metric
$\hat g$ is to have constant negative scalar curvature.  We shall not attempt
to give a survey of the results in either the positive case or the flat case.
Instead, we direct the reader to the introductions of the recent papers [MP],
[MPU], and the references contained therein.  We only note that a major
difference
between the negative case and the positive case (and the flat case) is the size
of $\Gamma$ for which there exists a Singular Yamabe metric.  For
example, suppose $\Gamma$ is a smooth submanifold of dimension $d$ in ${\bf
S}^n$. The results in [Av], [De], [LN], [Mz], [SY1], [SY2] then show there can
exist a complete metric with constant positive scalar curvature or zero scalar
curvature only if $d \leq {n-2  \over 2}$, and there can exist a complete metric
with constant scalar negative  curvature only if $d > {n-2 \over 2}$.

\newpar
The first results on the negative case appeared in Loewner and Nirenberg's
 seminal paper [LN] on partial differential equations invariant under conformal
transformations. In that paper, Loewner and Nirenberg showed when $\Gamma$ is a
smooth submanifold of dimension $d$ in ${\bf S}^n$ there exists a negative
Singular Yamabe metric $\hat g$ on ${\bf S}^n\backslash\Gamma$ if $d > {n-2
\over 2}$, and there does not exists such a metric if the Hausdorff dimension of
$\Gamma$ is less than ${n-2 \over 2}$.  Subsequent work by Aviles and McOwen in
[AM2] has since generalized the existence result in [LN] to an arbitrary compact
manifold proving that there exists a negative Singular Yamabe metric on
$M\backslash\Gamma$ when $\Gamma$ is a smooth submanifold of dimension $d$ if
and only if $d > {n-2 \over 2}$. This necessary and sufficient condition for
existence has been recently extended to smooth submanifolds of dimension $d$
with boundary (cf.\ [Fn1]).

\newpar
The original purpose of this work was to examine how singularities in the
structure of $\Gamma$ affect the behavior of a negative Singular Yamabe
metric  as a prelude to establishing existence.  But, in the course of our
investigations, it became apparent that singularities not only affect
behavior but also existence.   Specifically, we found that the dimension of the
{\it tangent cone} affects existence.

\result {Theorem}{Suppose that $\Gamma$ has a proper tangent cone at $p$, and
that the dimension of this tangent cone is less than ${n-2 \over 2}$.  Then
there can not exist a negative Singular Yamabe metric on $M\backslash\Gamma$.}

\noindent
We will define a proper tangent cone and discuss the dimension of a tangent
cone in the next section.  However, to  understand the significance of  our
result, the formal definitions of these terms are unnecessary.  For
instance, consider a singular hypersurface $\Gamma$ that is locally (near some
point) of the form
$$\{(x_1,\cdots,x_n): x_1^2 + \cdots + x_k^2 = x_n^3\}.$$
Our main result then shows that there can not exist a
negative Singular Yamabe metric on $M\backslash\Gamma$ when $k > {n+2 \over
2}$.  (The tangent cone to such a set at the origin is an upper half plane of
dimension $n-k$.)  This result is somewhat surprising for this choice of
$\Gamma$, since after all $\Gamma$ is a hypersurface and all the previous
results
pointed towards existence depending only on the ``global dimension'' of
$\Gamma$.

\newpar
We obtain our result by examining the behavior of the maximal positive solution
to the following semi-linear elliptic problem:
$$\left\{\eqalign{\Delta_gu &= u^q + Su \quad\hbox{on}\quad M\backslash\Gamma
\cr u(x) &\to +\infty \quad\hbox{as}\quad x \to \Gamma,}\right.\leqno(\dag)$$
where $\Delta_g$ is the Laplace operator on $(M,g)$, $q$ is an arbitrary
constant greater than one, and $S \in C^\infty(M)$. This works because
the equation in $(\dag)$ guarantees that the conformal metric $\hat g =
u^{q-1}\,g$ has constant negative scalar curvature (cf.\ [LP] or [Kz]), when
$q = {n+2 \over n-2}$ and $S$ is a specific multiple of the scalar curvature of
$g$.  Furthermore, if we can show that a solution $u$ to $(\dag)$ tends to
infinity at a sufficiently fast rate then the metric $\hat g$ will be complete.
Conversely, if we can show that the maximal solution to $(\dag)$ does not tend
to infinity at a fast rate then the metric $\hat g$ can not be complete, and
thus there can not exist a negative Singular Yamabe metric on
$M\backslash\Gamma$.

\vfill\eject

\line {\bf 2.\ Proper Tangent Cones \hfil}

\newpar
As stated previously, the original purpose of this paper was to study how
singularities in the structure of $\Gamma$ affect the behavior of the maximal
solution to $(\dag)$ and this the existence of a negative Singular Yamabe
metric.  We were going to start by considering the case where $\Gamma$ is a
stratified set, that is when $\Gamma$ can be decomposed into a finite number of
open embedded smooth submanifolds.  A prime example of such a set is an
algebraic variety of ${\bf R}^n$, or any set that can locally be obtained (in
some coordinate system) as an algebraic variety.  This case seemed to be the
natural first step after smooth submanifolds before a general closed subset.

\newpar
Through our investigation of this case, it became apparent that the behavior of
the maximal solution to $(\dag)$ depended on the tangent structure of $\Gamma$.
Therefore, instead of using the standard conditions on a stratified set (cf.\
[Wh1], [Wh2], [GM] or [Th]), we developed the following notion of a {\it proper
tangent cone} and placed conditions on the set through the admissibility of a
proper tangent structure in [Fn2].  As will become clear through the definitions
below, we do not require the set to be stratified for our main result.  However,
for the purpose of understanding the following definitions, it is useful to
restrict our attention to stratified sets,  partly because we require a proper
tangent cone to admit a stratification by smooth submanifolds.

\newpar
A stratified set is a set $\Gamma$ in an ambient manifold that can be
decomposed into a locally finite collection of disjoint smooth submanifolds.
This collection of smooth submanifolds is called a stratification of $\Gamma$,
and each submanifold is called a stratum (plural strata).  Typically, one also
enforces the condition that the stratification of $\Gamma$,
$\{\Sigma_\alpha\}$, satisfies the {\it axiom of frontier}: $\Sigma_\alpha \cap
{\rm clos}(\Sigma_\beta) \neq \emptyset$ implies that $\Sigma_\alpha \subset
{\rm clos}(\Sigma_\beta)$.

\newpar
Our first guess for defining a tangent cone to a stratified set $\Gamma$ at a
point $p$ is to generalize the definition of the tangent space to an embedded
smooth submanifold as the set of vectors $v \in T_pM$ such that there
exists a $C^1$ path $\gamma\colon[0,1] \to M$ with $\gamma(0) = p$, $\gamma'(0)
= v$ and $\gamma(t) \in \Gamma$ for all $t \in [0,1]$.  From this
definition it follows that if $v$ is a tangent vector to $\Gamma$ at $p$ then
for all $t \geq 0$ the vector $tv$ is also a tangent vector to $\Gamma$ at $p$.
Therefore, the set of tangent vectors to $\Gamma$ at $p$ is the cone over a set
of unit vectors in $T_pM$, which we call the {\it cone of tangent vectors}.
Notice that in this definition, if $p$ is a regular point of $\Gamma$, when the
set $\Gamma$ is locally a smooth submanifold near $p$, the cone of tangent
vectors coincides with the usual definition of the tangent space to $\Gamma$ at
$p$.

\newpar
We note however that for an arbitrary stratified set $\Gamma$ this notion of
the a cone of tangent vectors does not necessarily contain all the infinitesimal
information about $p$.  For example, consider the stratified set in ${\bf R}^2$
given by
$$\{(0,t): -1 \leq t \leq 1\}\ \bigcup\ \{(s,\sin(s)): s > 0\}.$$
Every tangent vector to this set at the origin is of the form $(0,y)$, but
every open set in $\Gamma$ containing the origin also contains points of the
form $(x,0)$.  To enlarge our notion of a tangent cone, we define the {\it
essential link of $\Gamma$ at $p$} as an infinitesimal version of an
$\epsilon$-link of $\Gamma$ at $p$.  In singularity theory, an $\epsilon$-link
of $\Gamma$ at $p$ is defined as the intersection of $\Gamma$ with a sphere of
radius $\epsilon$ centered at $p$, and is used to determine the topological
type of the singularity of $\Gamma$ at $p$ (cf.\ [Mi]).

\newpar
We define the essential link of $\Gamma$ at $p$ by examining the set of unit
vectors in $T_pM$ which are mapped into $\Gamma$ by geodesics.
Thus, our definition employs the map from $M\backslash\{p\}$ to $T_pM$ defined
by
$$\omega_p(x) = {v\over|v|}\quad\hbox{if}\quad x = \exp_p(v).$$
This map is well-defined in the punctured ball of radius $r_0$ about $p$,
$B_{r_0}(p)\backslash\{p\}$, where $r_0$ is a chosen constant less than the
injectivity radius of $(M,g)$.  Notice the map $\omega_p$ just gives the
spherical projection for a point $x$ in a geodesic polar coordinate system based
at $p$.  Using this map $\omega_p$, we can view an $\epsilon$-link of $\Gamma$
at $p$ as a subset of ${\bf S}^{n-1}$.  Therefore, we can take the limit of the
$\epsilon$-links of $\Gamma$ at $p$ as $\epsilon$ tends to zero.  The resulting
limit set could then be viewed as an infinitesimal version of an
$\epsilon$-link.  But to guarantee that this limits exists, it is more
convenient to define this {\it essential link of $\Gamma$ at $p$} as the
following subset of the unit sphere in $T_pM$:  $$L_p\Gamma = \bigcap_{0 <
\epsilon < r_0} {\rm clos} \big(\omega_p(\Gamma\cap B_\epsilon^*)\big)$$
where $B_\epsilon^*$ is the punctured ball of radius $\epsilon$,
$B_\epsilon(p)\backslash\{p\}$.

\newpar
This definition of an essential link is closely related to the definition
of an tangent cone in geometric measure theory (cf.\ [Fe] or [Mo]).  In fact,
in geometric measure theory, the tangent cone of an arbitrary set is defined to
be the cone over the essential link.  This is the definition of a tangent cone
we will be using in this paper.  Specifically, we define the {\it tangent
cone of
$\Gamma$ at $p$} as
$$T_p\Gamma = \{w \in T_pM\colon w = tv \hbox{ with } t \geq 0 \hbox{ and } v
\in L_p\Gamma\}.$$
We note that this cone of tangent vectors to $\Gamma$ is clearly contained in
$T_p\Gamma$, and for most nice stratified sets the cone of tangent vectors is
equal to $T_p\Gamma$.

\newpar
In defining a {\it proper tangent cone}, we would like to guarantee that both
$T_p\Gamma$ and $L_p\Gamma$ satisfy certain regularity conditions.  This is
because it is possible to construct stratified sets where $L_p\Gamma$ (and
hence $T_p\Gamma$) is a Cantor set.  Simply notice that the process for
defining the essential link is similar to an inductive process.  Therefore, it
is not hard to construct a stratified set (with infinite topological type) that
has an essential link that is a Cantor set.  To avoid such sets, we demand that
the essential link is a stratified set when defining a proper tangent cone.

\state {{\bf Definition:} A tangent cone to $\Gamma$ at a point $p$
is called a {\it proper tangent cone} if the essential link of $\Gamma$ at $p$
admits a stratification by smooth submanifolds.}

\noindent
By requiring the essential link to be stratified, we guarantee that a proper
tangent cone is also stratified.

\newpar
We define the dimension of a stratified set $\Gamma$ to be $d$ is $d$
is the dimension of the largest submanifold in a stratification of $\Gamma$.
Therefore, we can define the dimension of a proper tangent cone and the
dimension of the essential link of a proper tangent cone using this definition
of the dimension of a stratified set.  From this, it follows that
$\dim(T_p\Gamma) = \dim(L_p\Gamma) + 1$, whenever $T_p\Gamma$ is a proper
tangent cone.  We note this concept of dimension does not necessarily agree
with the Hausdorff dimension of $T_p\Gamma$ as a subset of $T_pM$.  For
example, the Hausdorff dimension of the set
$$\Gamma = \{(0,t): -1 \leq t \leq 1\}\ \bigcup\ \{(s,\sin(s)): s > 0\}\subset
{\bf R}^2$$
is two, while the dimension we are using for $\Gamma$ as a stratified set says
the dimension is one.

\vfill\eject

\line {\bf 3.\ Outline of Proof \hfil}

\newpar
Let $\beta_0 = {2 \over q-1}$ and $d_0 = n-2 -\beta_0$.  Notice when $q = {n+2
\over n-2}$, we have $\beta_0 = d_0 = {n-2 \over 2}$.  We will first show that
if the maximal positive solution $u$ to $(\dag)$ has a strong singularity at
$\Gamma$, that is $u$ satisfies the estimate
$$0 < C_1 \leq \rho(x)^{\beta_0}\,u(x) \leq C_2 < +\infty\leqno(3.1)$$
where $\rho(x) = dist_g(x,\Gamma)$ then the asymptotic behavior of $u$ near a
point $p \in \Gamma$ can be described in geodesic polar coordinates $(r,\omega)$
at $p$ by
$$v(\omega) = \lim\limits_{r \to 0}\ r^{\beta_0}\,u(r,\omega)\leqno(3.2)$$
with $v$ being the maximal positive solution to
$$\left\{\eqalign{\Delta_\omega v = v^q + \beta_0\,d_0\,v \quad\hbox{on}\quad
{\bf S}^{n-1}\backslash L_p\Gamma \cr v(\omega) \to + \infty
\quad\hbox{as}\quad \omega \to L_p\Gamma,}\right.\leqno(\ast)$$
where $\Delta_\omega$ represents the Laplace operator on ${\bf S}^{n-1}$
with respect to the standard metric. The fact that $u$ has a strong
singularity at $\Gamma$ implies that $v$ must also have a strong singularity at
$L_p\Gamma$; $v$ must satisfy
$$0 < C_1 \leq \sigma(\omega)^{\beta_0}\,v(\omega) \leq C_2 < +\infty
\quad\hbox{as}\leqno(3.3)$$
where $\sigma(\omega) = dist(\omega,L_p\Gamma)$ on ${\bf S}^{n-1}$.  Such a
solution to $(\ast)$ can only exist if the dimension of $L_p\Gamma$ is greater
than $d_0-1$ (see Lemma 2 in the next section). Hence, the dimension of the
tangent cone $T_p\Gamma$ must be greater than $d_0$ whenever the maximal
solution
to $(\dag)$ has a strong singularity at $\Gamma$.

\newpar
We note that a solution to $(\dag)$ with a strong singularity at $\Gamma$
implies that the conformal metric $\hat g = u^{q-1}\,g$ is a complete metric.
Therefore, we conclude from the above analysis that any obstruction to the
existence of a negative Singular Yamabe metric occurs when the dimension of
$T_p\Gamma$ is less than or equal to ${n-2 \over 2}$. We would like to conclude
that if the maximal solution does not have a strong singularity then the metric
$\hat g$ is not complete.  But, our examination of the asymptotic behavior of
solutions to $(\dag)$ does not allow this conclusion.  Our analysis is based
solely on how the dimension of $T_p\Gamma$ affects the behavior.

\newpar
We obtain our main result follows from a delicate asymptotic analysis which
shows that when the dimension of $T_p\Gamma$ is less than $d_0$ there exists a
$\delta > 0$ and an open path $\gamma\colon (0,1) \to M\backslash\Gamma$ with $p
= \lim\limits_{t \to 0} \gamma(t)$ such that the maximal solution to $(\dag)$
satisfies
$$\limsup\limits_{x \in \gamma} \rho(x)^{\beta_0 - \delta}\,u(x) < +
\infty.\leqno(3.4)$$
The last step in this analysis is to compare the asymptotic behavior of
the maximal solution to $(\dag)$ to a solution of
$$\left\{\eqalign{&\Delta_\omega v = \beta_0d_0\,v \quad\hbox{on}\quad {\bf
S}^{n-1}\backslash\L_p\Gamma \cr &v(\omega) \to +\infty \quad\hbox{as}\quad
\omega \to L_p\Gamma.}\right.\leqno(\ast)$$
This is the step that requires the dimension of $T_p\Gamma$ to be strictly less
than $d_0$. Once we have established (3.4), a simple calculation shows that the
metric $\hat g = u^{q-1}\,g$ is not  complete.  Hence, when the dimension of a
proper tangent cone is less than ${n-2 \over 2}$ there can not exist a negative
Singular Yamabe metric.

\vfill\eject
\line {\bf 4.\ Analytical Preliminaries \hfil}

\newpar
To obtain our main result, we need the following analytic results on the
behavior of positive solutions to $\Delta_g u = u^q + Su$  proved in [Av], [FM]
and [Fn1] (see also [Fn2]).  The proofs of these results all rely on standard
elliptic theory (cf.\ [GT]).  We will also need to use the method of upper and
lower solutions, and a global Green's function (cf.\ [Au], [LP], and [Mc2]).

\newpar
The first important result on the behavior of positive solutions
to $(\dag)$ is the following {\it apriori} upper bound found in [Av] and [Vn]

\result {Lemma 1}{Suppose $u$ is a positive solution of $\Delta_gu = u^q + Su$
in $M\backslash\Gamma$. Then
$$\rho(x)^{\beta_0}\,u(x) \leq C$$
for some positive constant $C$ independent of $u$ and $\Gamma$.}

\noindent
In both [Av] and [Vn], this result is stated for the case where $\Gamma$ is
a smooth submanifold, but the proofs are valid for an arbitrary closed set.
One of the first uses of this apriori upper bound is the following nonexistence
result proved in [Av] and [Vn], and extended in [Fn1], [Fn2] to the following
form:

\result {Lemma 2}{Suppose $p$ is a regular point of $\Gamma$ and the
dimension of $\Gamma$ at $p$ is less than or equal to $d_0$.  Then there exists
no positive solution of $\Delta_gu = u^q + Su$ in $M\backslash\Gamma$ satisfying
$u(x) \to \infty$ as $x \to p$.}

\noindent
This result was originally stated for smooth submanifolds, but by introducing a
cut-off function on $\Gamma$ around $p$ into the Aviles' proof one obtains the
above generalization (see [Av], [Fn1], [Fn2]).  The regularity on $\Gamma$ is
needed to establish a nice coordinate system about $\Gamma$, so one may estimate
integrals near $\Gamma$.

\newpar
We obtain more refined results on the asymptotic behavior of the maximal
solution to $(\dag)$ from a {\it Harnack-type inequality} and {\it
derivative estimates} for $\Delta_gu = u^q + Su$ (cf.\ [FM] and [Fn2]). These
estimates are obtained directly from standard elliptic theory as a consequence
of the apriori upper bound.  Hence, these results do not depend on the
regularity of the $\Gamma$.

\result {Lemma 3}{Let $u$ be a positive solution to $\Delta_g u = u^q + Su$ in
$M\backslash\Gamma$ and $x_0$ be a fixed point in $M\backslash\Gamma$.  Then
there exists a positive constant $C_1$ independent of $u$, $x_0$ and $\Gamma$
such that
$$\sup\limits_{y \in B}u(y) \leq C_1\,\inf\limits_{y \in B} u(y)$$
where $B = \{y \in M\backslash\Gamma: dist_g(x_0,y) \leq {1\over
8}\rho(x_0)\}$.}

\noindent
To state the derivative estimates, it will be convenient to
introduce a weighted norm on functions in $C^{2,\alpha}(\Omega)$  where
$\Omega$ is a smooth submanifold of dimension $n$ with (smooth) boundary in
$M$.  Let $d_\xi = dist_g(\xi,\partial\Omega)$ and $d_{\xi,\eta} =
\min(d_\xi,d_\eta)$. Furthermore, let ${\cal D}u(\xi)$ and ${\cal D}^2u(\xi)$
represent the covariant derivative of $u$ and ${\cal D}u$ with respect to an
orthonormal frame at $\xi$.  We then define the norm $|\cdot|^\ast_{2,\alpha}$
of a function $u \in C^{2,\alpha}(\Omega)$ by
$$\left|u\right|^\ast_{2,\alpha;\Omega} = \left|u\right|^\ast_{2,0;\Omega} +
\left[u\right]^\ast_{2,\alpha;\Omega}\leqno(4.1)$$
where
$$\left|u\right|^\ast_{2,0;\Omega} = \sup\limits_{\xi\in\Omega} |u(\xi)| +
\sup\limits_{\xi\in\Omega} d_\xi\,|{\cal D}u(\xi)| + \sup\limits_{\xi\in\Omega}
d_\xi^2\,|{\cal D}^2u(\xi)|\leqno(4.2)$$
and
$$\left[u\right]^\ast_{2,\alpha;\Omega} = \sup\limits_{\xi,\eta\in \Omega}
d_{\xi,\eta}^{2+\alpha} \,{|{\cal D}^2u(\xi) - {\cal D}^2u(\eta)| \over
|\xi-\eta|_g^\alpha}.\leqno(4.3)$$
In (4.2), $|{\cal D}u(\xi)|$ is the length of the covariant
derivative of $u$ at $\xi$, and $|{\cal D}^2u(\xi)|$ is the length of the
covariant derivative of ${\cal D}u$ at $\xi$. In (4.3), the addition of vectors
in different tangent spaces is accomplished by using parallel transport.

\result {Lemma 4}{Let $u$ be a positive solution to $\Delta_g u = u^q + Su$ on
$M\backslash\Gamma$ and let $x_0$ be a fixed point in $M\backslash\Gamma$.  Then
there exists a positive constant $C_2$ independent of $x_0$ and $\Gamma$ such
that
$$|u|^\ast_{0,2;B} \leq C_2 |u|_B$$
where $B = \{y \in M\backslash\Gamma: dist_g(x_0,y) \leq {1\over
8}\rho(x_0)\}$.}

\noindent
From these derivative estimates, we  derived the uniqueness of a
positive solution to $(\dag)$ with a strong singularity at $\Gamma$ in [Fn1].

\result {Lemma 5}{There exists at most one positive solution to $(\dag)$ with a
strong singularity at $\Gamma$.}

\noindent
The proof of this uniqueness also requires the Asymptotic Maximum Principle of
Cheng and Yau (cf.\ [CY] and [Au]).

\newpar
Probably, the most unusual aspect of our nonexistence proof is the use of
existence methods.  Specifically, we make use of the method of upper and lower
solutions and the existence of a global Green's function.  We recall that in
the method of upper and lower solutions one looks for an upper solution
$\varphi$, a function that satisfies
$$L[\varphi] + f(x,\varphi) \leq 0,$$
and a lower solution $\psi$, a function that satisfies
$$L[\psi] + f(x,\psi) \geq 0$$,
 such that $\varphi \geq \psi$.  A monotone iteration scheme then  guarantees
the existence of a solution $u$ of $L[u] + f(x,u) = 0$ which satisfies $\psi
\leq u \leq \varphi$.  However, in our proof, it is necessary to work with
Holder
continuous upper and lower solutions, instead of the usual $C^2$ functions.
Therefore, we use Calabi's extension of the Hopf maximum principle in the method
of upper and lower solutions (cf.\ [Ca]).  This extension is based on the
following definition for a continuous function to weakly satisfy a differential
inequality.

\state {{\bf Definition:} Let $L$ be a uniformly elliptic operator of
divergence type.  We say a continuous function $u$ in an open domain $\Omega$
satisfies
$$L[u] \geq v \quad\hbox{weakly in any subset $K$ of $\Omega$}$$
if for each point $x_0 \in K$ and any given positive constant $\epsilon$ there
exists an open neighborhood $V_{x_0,\epsilon} \subset \Omega$ containing $x_0$
and a $C^2$ function $u_{x_0,\epsilon}$ such that the difference function $u-
u_{x_0,\epsilon}$ achieves it minimum value at $x_0$ and $u_{x_0,\epsilon}$
satisfies the inequality
$$L[u_{x_0,\epsilon}] \geq v-\epsilon \quad\hbox{in}\quad V_{x_0,\epsilon}$$
in the usual sense. Similarly, we say that $L[u] \leq v$ weakly if $L[-u] \geq
-v$ weakly, as defined above.}

\noindent
In fact from the methods in [AM1], we can reduce the method of upper and lower
solutions for $(\dag)$ to a search for nonnegative weak lower solution $w$ that
satisfies $w(x) \to +\infty$ as $x \to \Gamma$. (see also [Fn2])

\newpar
As noted in the outline of our proof, we compare the asymptotic behavior of the
maximal solution  to $(\dag)$ to a solution of the linear
equation$\Delta_\omega v = \beta_0d_0v$ on ${\bf S}^{n-1}\backslash
L_p\Gamma$.  We construct solutions to this linear equation by using the global
Green's function to $-\Delta_\omega + \beta_0d_0$ on ${\bf S}^{n-1}$ to define a
Poisson transform (cf. [Fn1]).  Specifically, there exists a global Green's
function ${\cal G}(x,y)$ when $d_0 > 0$ and this Green's function satisfies the
estimates $$C_1\,|x-y|_\omega^{3-n} \leq {\cal G}(x,y) \leq
C_2\,|x-y|_\omega^{3-n}$$ when $n\geq 4$, and
$$C_1\,\big(1 + \big|\log|x-y|_\omega\big|\big) \leq {\cal G}(x,y) \leq C_2\,
\big(1 + \big|\log|x-y|_\omega\big|\big)$$
when $n=3$, where $C_1,C_2$ are positive constants (depending on $n$ and $q$),
and $|x-y|_\omega$ is the distance between $x$ and $y$ on ${\bf S}^{n-1}$.
Recall that for every $y in {\bf S}^{n-1}$, the function $u(x) = {\cal
G}(x,y)$ is the positive distributional solution of $-\Delta_\omega v +
\beta_0d_0 v = \delta_y$ with $\delta_y$ is the Dirac delta distribution
at $y$. With this Green's function, we define a Poisson transform on ${\bf
S}^{n-1}\backslash L_p\Gamma$ by
$$T[\varphi](x) = \int\limits_{L_p\Gamma} {\cal
G}(x,y)\,\varphi(y)\,d\nu(y)$$
where $\nu(y)$ is a measure on $L_p\Gamma$ compatible with the induced
Riemannian measure on each stratum.  Notice for Lebesgue integrable functions
(with respect to the measure $\nu$) the Poisson transform constructs classical
solutions of $\Delta_\omega v = \beta_0d_0v$ on ${\bf S}^{n-1}\backslash
L_p\Gamma$.

\vfill\eject
\line {\bf 5.\ Asymptotic Behavior of a Strong Singularity \hfil}

\newpar
In this section, we complete the first half of the proof of our main result by
examining the asymptotic behavior of a solution to $\dag)$ with a strong
singularity.  In particular, we prove

\result {Theorem A}{Suppose that $\Gamma$ has a proper tangent cone at $p$, and
let $u$ be the unique positive solution to $(\dag)$ with a strong singularity
at $\Gamma$.  Then the dimension of $T_p\Gamma$ must be greater than $d_0$.
Moreover, the function $v$ defined by $(3.2)$ is the unique positive solution
to $(\ast)$ with a strong singularity at $L_p\Gamma$.}

\noindent
This theorem describes the unique solution to $(\dag)$ with a strong
singularity at $\Gamma$ in terms of another unique solution with a strong
singularity.  Therefore, we may use this theorem inductively, to prove that
very strong geometric conditions are needed to prove the existence of a
solution to $(\dag)$ with a strong singularity at $\Gamma$, provided the
essential link $L_p\Gamma$ satisfies sufficient regularity conditions.
Consequently, any existence proof for a negative Singular Yamabe metric becomes
much more complicated; the previous existence proofs are all based on proving
the existence of a positive solution to $\dag)$ with a strong singularity at
$\Gamma$, since a solution with a strong singularity guarantees that the
conformal metric $\hat g = u^{q-1}\,g$ is a complete metric.

\newpar
We will first establish some notation before proving Theorem A.  Let
$(r,\omega)$ be a geodesic coordinate system centered at $p$. In terms of this
coordinate system, the operator $\Delta_g$ takes on the form:
$$\Delta_g = \partial_r^2 + {n-1 \over
r}\partial_r + {1 \over r^2}\Delta_\omega + {\cal L}\leqno(5.1)$$
where ${\cal L}$ is a second order perturbation operator of form
$${\cal L} = {\cal O}(r)\,\partial_r + {\cal O}(r^{-1})\,\partial_\omega^2 +
{\cal O}(r^{-1})\,\partial_\omega.\leqno(5.2)$$
In (5.2), we use ${\cal O}(r^k)$ to mean the term is dominated by $C\,r^k$ as
$r \to 0$, i.e.\ big-oh notation.  Since we will constantly be using
expressions similar to (3.1) and (3.3), we introduce the following notation
$$f(r) \approx g(r) \quad\hbox{as}\quad r \to 0$$
to mean that there exist fixed positive constants $C_1,C_2$ such that
$$C_1\,f(r) \leq g(r) \leq C_2\,f(r) \quad\hbox{as}\quad r \to 0.$$
Thus, a solution to $(\dag)$ with a strong singularity can be described as
$u(x) \approx \rho(x)^{-\beta_0}$.

\newpar
The proof of Theorem A is divided into three propositions concerning the
behavior of the following two functions
$$\left\{\eqalign{&\overline{v}(\omega) = \limsup\limits_{r\to 0}\
r^{\beta_0}\,u(r,\omega) \cr &\underline{v}(\omega) = \liminf\limits_{r \to 0}\
r^{\beta_0}\,u(r,\omega)}\right.$$
which will show that $\overline{v}\equiv\underline{v}$, and thus that the
function defined by (3.2) exists.  These propositions together prove that $v$
is the unique positive solution of $(\ast)$ with a strong singularity at
$L_p\Gamma$, from which Theorem A follows by appealing to Lemma 2 in the
previous section.

\result {Proposition 1}{$\overline v,\ \underline{v} \in
C^{0,\alpha}({\bf S}^{n-1}\backslash L_p\Gamma)$ for some $\alpha \in
(0,1)$.}

\noindent
{\bf Proof:}  Choose $\Omega \subset\subset {\bf S}^{n-1}\backslash
L_p\Gamma$, and define $\Omega(r) = \{x \in M: x=\exp_p(rv) \hbox{ where } v
\in \Omega\}$. Then there exists  $\delta>0$ and
$\epsilon>0$ sufficiently small such that $dist_g(\Omega(r),\Gamma) \geq
\epsilon r$ for all $r \in (0,\delta)$,  by definition of $L_p\Gamma$.
From the apriori upper bound in Lemma 1 there exists a positive constant $C$
depending on $\epsilon$ and $\delta$ such that
$$v(r,\omega) = r^{\beta_0}\,u(r,\omega) \leq C
\quad\hbox{when}\quad \omega  \in \Omega$$ for all $r \in (0,\delta)$.  The
derivative estimates in Lemma 4 then imply $v_r(\omega) = v(r,\omega) \in
C^{2,\alpha}(\overline{\Omega})$ for all $r \in (0,\delta)$ and there exists a
positive constant $C$ such that $||v_r||_{2,\alpha;\Omega} \leq C$, (i.e.\
$||v_r||_{2,\alpha;\Omega}$ is bounded independent of $r$), where
$||\cdot||_{2,\alpha;\Omega}$ is the standard norm on the function space
$C^{2,\alpha}(\overline{\Omega})$.  Since the derivatives in $\Omega$ are
bounded independent of $r$, it follows that the H\"older semi-norm
$[v]_{0,\alpha;\Omega}$ is bounded independent of $r$.  Therefore, it follows by
standard arguments that ${\overline v}, {\underline v} \in
C^{0,\alpha}(\overline{\Omega})$, and thus ${\overline v}, {\underline v} \in
C^{0,\alpha}({\bf S}^{n-1}\backslash L_p\Gamma)$.  \qed

\newpar
\result {Proposition 2}{$\overline{v}$ satisfies $\Delta_\omega v \geq v^q +
\beta_0d_0v$ weakly on ${\bf S}^{n-1}\backslash L_p\Gamma$, and $\underline{v}$
satisfies $\Delta_\omega v \leq v^q + \beta_0d_0v$ weakly on ${\bf
S}^{n-1}\backslash L_p\Gamma$.}

\noindent
{\bf Proof:} This proposition is proved by appealing to a one-dimensional version
of the asymptotic maximum principle, and the definitions of $\overline v$ and
$\underline v$.  First, we note that under the change of variables
$r=e^{-t}$ and $v(t,\omega) = e^{\beta_0t}\,u(t,\omega)$ the equation $\Delta_gu
= u^q + Su$ may written as
$$\partial_t^2v + (\beta_0-d_0)\partial_tv + \Delta_\omega v = v^q +
\beta_0d_0 v + {\cal L}v$$
where ${\cal L}v = {\cal O}(e^{-t}) D^2v + {\cal O}(e^{-t})Dv + {\cal
O}(e^{-2t})v$.  With this change of variables, we find
$$\left\{\eqalign{\overline{v}(\omega) = &\limsup\limits_{t \to \infty}
v(t,\omega) \cr \underline{v}(\omega) = &\liminf\limits_{t \to \infty}
v(t,\omega).}\right.$$
Next, we choose $\Omega \subset\subset {\bf S}^{n-1}\backslash L_p\Gamma$ as in
Proposition 1, so there exists $T$ sufficiently large such that $v_t(\omega) =
v(t,\omega) \in C^{2,\omega}(\overline{\Omega})$ for all $t \in [T,\infty)$ and
there exists a positive constant $C$ such that $||v_t||_{2,\alpha} \leq C$
uniformly for all  $t \in [T,\infty)$.  From a one-dimension version of the
asymptotic maximum principle and the definition of $\overline v$, we can
choose a monotonically increasing sequence $\{T_n\}$ with $T_n \to  \infty$
such that
$$\left\{\eqalign{&\overline{v}(\omega_0) = \lim\limits_{n\to
\infty} v(T_n,\omega_0) \cr &\lim\limits_{n\to\infty} \partial_t
v(T_n,\omega_0) = 0 \cr &\lim\limits_{n\to\infty} \partial_t^2
v(T_n,\omega_0) \leq 0,}\right.$$
for any point $\omega_0 \in \Omega$. The Arzela-Ascoli theorem (and the
compactness of the embedding $C^{2,\alpha} \to C^2$) implies there exists a
subsequence $\{v_m(\cdot) = v(T_{n_m},\cdot)\}$ converging to $v \in
C^2(\omega)$.  We have $v \leq \overline{v}$ and $v(\omega_0) =
\overline{v}(\omega)$.  Therefore, for any $\epsilon > 0$, we can find a
neighborhood $U$ of $\omega_0$ such that
$$\Delta_\omega v \geq \overline{v}^q +
\beta_0d_0\,\overline{v}-\epsilon \quad \hbox{in} \quad U.$$
The function $\overline{v} \in C^{0,\alpha}({\bf S}^{n-1}\backslash
L_p\Gamma)$ then is a weak lower solution of $\Delta_\omega v = v^q +
\beta_0d_0 v$ as desired.

\newpar
The proof that $\underline v$ is a weak lower solution to $\Delta_\omega v =
v^q + \beta_0d_0v$ is based on repeating the above argument, using that
$$\left\{\eqalign{&\underline{v}(\omega_0) = \lim\limits_{n\to
\infty} v(T_n,\omega_0) \cr  &\lim\limits_{n\to\infty} \partial_t
v(T_n,\omega_0) = 0 \cr &\lim\limits_{n\to\infty} \partial_t^2
v(T_n,\omega_0) \geq 0}\right.$$
for some monotonically increasing sequence $\{T_n\}$ with $T_n \to \infty$. \qed

\result {Proposition 3}{$\overline{v}(\omega) \approx \underline{v}(\omega)
\approx \sigma(\omega)^{-\beta_0}$ as $\sigma(\omega) \to 0$.}

\noindent
{\bf Proof:} To prove this proposition, choose $\epsilon > 0$ and let
$N_\epsilon = \{v \in {\bf S}^{n-1}: \sigma(v) < \epsilon\}$,
$U_\epsilon = \{v \in {\bf S}^{n-1}: \sigma(v) > 2\epsilon\}$. Then there exists
$\delta > 0$ depending only on $\epsilon$ such that $\Gamma \cap B_p(\delta) =
\{x \in \Gamma: dist_g(x,p) < \delta\}$ is contained in $N_\epsilon(\delta) =
\{x = \exp_p(rv): 0 \leq r < \delta,\ v \in N_\epsilon\}$. It then follows
for $x \in U_\epsilon$ that there exists positive constants $C_1,C_2$
independent of $\epsilon$ and $\delta$ with
$$C_1\,r\,\sigma(v) \leq dist_g(\exp_p(rv),N_\epsilon(\delta))\leq
C_2\,r\,\sigma(v)$$
and
$$C_1\,\rho(\exp_p(rv)) \leq dist_g(\exp_p(rv),N_\epsilon(\delta)) \leq
C_2\,\rho(\exp_p(rv)).$$
Thus, it follows that $r\sigma(v) \approx \rho(\exp_p(rv)$ as $r \to 0$, and
therefore that
$$r^{\beta_0}\,u(r,\omega) \approx \sigma(\omega)^{-\beta_0}$$
for $\omega \in \Omega$. Passing through the limits, we see that
$$\overline{v}(\omega) \approx \underline{v}(\omega) \approx
\sigma(\omega)^{-\beta_0}$$
which complete the proof of the proposition. \qed

\noindent
{\bf Proof of Theorem A:} By definition of $f \approx g$ and ${\overline
v},{\underline v}$, there exists a small positive constant $c$ such that
$$c\overline{v}(\omega) \leq \underline{v}(\omega) \leq \overline{v}(\omega)
\quad\hbox{for all}\quad \omega \in {\bf S}^{n-1}\backslash L_p\Gamma.$$
We also have $c\overline{v}$ to be lower solution to $\Delta_\omega v = v^q +
\beta_0d_0v$, so the method of upper and lower solutions guarantees the
existence of a positive solution to $\Delta_gv = v^q + \beta_0d_0v$ with a
strong singularity at $L_p\Gamma$ that satisfies $c\overline{v} \leq v \leq
\underline{v}$, and the existence of a positive solution to
$\Delta_\omega v = v^q  + \beta_0d_0v$ with a strong singularity at $L_p\Gamma$
with $v \geq \overline{v}$.  The uniqueness of solutions with a strong
singularity (Lemma 5) then shows that $\overline{v} \equiv \underline{v}$, and
thus the limit
$$v(\omega) = \lim\limits_{r \to 0} \ r^{\beta_0}\,u(r,\omega)$$
exists and is the unique positive solution of $\Delta_\omega v = v^q +
\beta_0d_0v$ in ${\bf S}^{n-1} \backslash L_p\Gamma$ with a strong singularity
at $L_p\Gamma$. The localization of Aviles' nonexistence result then shows that
the local dimension of $L_p\Gamma$ must be greater than $d_0-1$, and thus the
dimension of the tangent cone must be greater than $d_0$. \qed

\newpar
Let us take a moment to discuss the existence of a solution to $(\dag)$ with a
strong singularity.  This is discussed in detail in [Fn2], and the outline for
the procedure is contained in [Fn1].  For existence, we demand that the cone of
tangent vectors of $\Gamma$ at $p$ is equal to the tangent cone of $\Gamma$ at
$p$ for all $p$.  Further, we demand that the tangent cones are compatible for
two points in the same stratification, and lastly we demand that we can replace
a portion of $\Gamma$ near $p$ with the tangent cone of $\Gamma$ at $p$.  A
modification of the methods in [Fn1] allows us to construct a solution to
$\Delta_gu=u^q + Su$ that is singular at the high dimensional portion of
$T_p\Gamma$ (for strata in $T_p\Gamma$ with dimension greater than $d_0$).  To
conclude that $u(x) \to +\infty$ as $x \to p$,  we further require that a
high-dimensional piece of $T_p\Gamma$ is cone-like, that is of the form
$$\{(x_1,\cdots,x_m,0,\cdots,0): x_1^2+\cdots+x_k^2 = x_m^2,\ x_m \geq 0\}
\subset {\bf R}^n$$
where $k < m \leq n$, with $m>d_0$.  This last condition is necessary so we may
apply the {\it continuation argument} developed in [Fn1].

\vfill\eject

\line {\bf 6.\ Proof of Main Result\hfil}

\newpar
The result of the previous section shows that the dimension of every proper
tangent cone must be greater than $d-0$ for the existence of a positive
solution to $(\dag)$ with a strong singularity at $\Gamma$. We would like to
claim that if the maximal solution $u$ to $(\dag)$ does not have a strong
singularity then the conformal metric $u^{q-1}\,g$ is not complete.  But, this
is difficult to prove.  Instead, we will show that if the dimension of
$T_p\Gamma$ is less than $d_0$, then the metric $\hat g$ is not complete.  The
critical case where the dimension of a proper tangent cone is equal to $d_0$
requires further study; our proof breaks down exactly in this case.

\result {Theorem B}{Suppose $\Gamma$  has a proper tangent cone at $p$ with
the dimension of this tangent cone less than $d_0$, and let $u$ be the maximal
positive solution to $(\dag)$.  Then the conformal metric $\hat g =
u^{q-1}\,g$ is not a complete metric on $M\backslash\Gamma$.}

\noindent
{\bf Proof:} It suffices to show when $d_0 > 0$ that there exists an $\epsilon
> 0$ such that every positive solution to $(\dag)$ satisfies
$$r^{\beta_0}\,u(r,\omega) = {\cal O}(r^\epsilon) \quad \hbox{as} \quad r
\to 0\leqno(6.1)$$
for all $\omega \in {\bf S}^{n-1}\backslash L_p\Gamma$.  Theorem A shows that
$$r^{\beta_0}u(r,\omega) \to 0 \quad \hbox{as}\quad r \to 0,\leqno(6.2)$$
since $\dim(T_p\Gamma) < d_0$.  To improve (6.2) to (6.1), we
first show that there exists a positive function $f$ such that
$$r^{\beta_0}\,u(r,\omega) \approx f(r) \quad \hbox{as}\quad r \to
0,\leqno(6.3)$$  with $f$ satisfying $f(r) \to 0$ as  $r \to 0$, and
$$\lim\limits_{r \to 0} {rf'(r) \over f(r)}
\quad\hbox{and}\quad \lim\limits_{r\to 0} {r^2f''(r) \over f(r)} \quad
\hbox{both exist.}$$
As a result of these conditions, we find that
$$\overline{v}(\omega) = \limsup\limits_{r \to 0} {r^{\beta_0}\,u(r,\omega)
\over f(r)}\leqno(6.4)$$
is a positive function on ${\bf S}^{n-1}\backslash L_p\Gamma$.  A simple
modification of the arguments in the proof of Propositions 1  in Section 5
shows that $\overline v \in C^{0,\alpha}({\bf S}^{n-1}\backslash
L_p\Gamma)$.

\newpar
We prove the existence of such a function $f$ as follows.  First choose
$\omega_0 \in {\bf S}^{n-1}\backslash L_p\Gamma$,  and define
$$f_0(t) = \sup\limits_{r \in [0,t]} \,r^{\beta_0}\,u(r,\omega_0).$$
Notice, we have defined $f_0$ so that $(6.3)$ is satisfied when $\omega =
\omega_0$. The Harnack inequality, Lemma 3 in Section 4, then guarantees that
(6.3) is satisfies for all $\omega in {\bf S}^{n-1}\backslash L_p\Gamma$.
Notice we have  defined $f_0$ so that it is Lipschitz continuous, and hence is
differentiable almost everywhere.  Choose a $C^1$ function $f_1$ close to $f_0$
in the sup-norm.  By computing the derivative of $f_0$ when possible, we find
that either
$$f_0'(r) = \partial_r (r^{\beta_0}\,u(r,\omega_0))$$
or $f_0'(r) = 0$, so the derivative estimates in Lemma 4 imply that
$$0 \leq {r\,f_0'(r) \over f_0(r)} < C$$
when $f_0'(r)$ is defined.  Therefore, we can arrange that
$$0 < {r\,f_1'(r) \over f_1(r)} < C.\leqno(6.5)$$
by choosing $f_1$ to be strictly increasing.  We may change the $C^1$
approximation to be either concave up or concave down without altering the truth
of (6.3) or (6.5). Thus, we can choose $f_1$ such that
$$\lim\limits_{r \to 0} {r\,f_1'(r) \over f_1(r)} \quad
\hbox{exists}.\leqno(6.6)$$
Further, we may choose $f_1$ such that $rf_1'(r) / f_1(r)$ is either
monotonically increasing or monotonically decreasing without altering the truth
of (6.3), (6.5), or (6.6).  Thus, by standard results (cf.\ [R]), we find
that $rf_1'(r) / f_1(r)$ is differentiable almost everywhere.  Repeating the
arguments in proceeding from $f_0$ to $f_1$, we may thus choose a $C^2$ function
$f$ satisfying the (6.3), (6.5), and (6.6) with
$$\left|{r^2 f''(r) \over f(r)}\right| < C.\leqno(6.7)$$
We may further choose $f''(r)$ to be monotone  without altering the truth of
(6.3), (6.5), (6.6) or (6.7) to ensure that
$$\lim\limits_{r \to 0} {r^2 f''(r) \over f(r)} \quad
\hbox{exists}.\leqno(6.8)$$
Therefore,  we have the existence of the desired function $f$.

\newpar
We now consider two cases in proving that (6.2) may be improved to (6.1):
$${\rm (i)}\ \lim\limits_{r\to 0} {rf'(r) \over f(r)} = \epsilon > 0
\quad\hbox{and}\quad {\rm (ii)}\ \lim\limits_{r\to 0} {rf'(r) \over f(r)} =
0.$$
In the first case, a comparison theorem implies that $f(r) \leq
C\,r^{\epsilon}$, and thus $u(r,\omega) \leq Cr^{-\beta_0+\epsilon}$ which
implies that $u^{q-1} \leq r^{-2+\epsilon'}$ and hence noncompleteness of the
metric $u^{q-1}\,g$.  In the second case, we find that
$$\lim\limits_{r \to 0} {r^2 f''(r) \over f(r)} = 0,$$
since if we assume
$$\lim\limits_{r \to 0} {r^2f''(r) \over f(r)} = \epsilon \neq 0
\quad\hbox{and}\quad \lim\limits_{r \to 0} {rf'(r) \over f(r)} = 0$$
then
$$\lim\limits_{r \to 0} r{d\over dr}\left({rf'(r) \over f(r)}\right) =
\epsilon$$
which implies
$$\left|{rf'(r) \over f(r)}\right| \approx -\log(r).$$
Thus ${rf'(r) \over f(r)}$ tends to infinity as $r \to 0$ contradicting
${rf'(r) \over f(r)} \to 0$.  Define the function $v(r,\omega)$ by  setting
$$u(r,\omega)  = r^{-\beta_0}\,f(r)\,v(r,\omega)$$
Computing $\Delta_gu$ using (5.1), we find that
$$\eqalign{\Delta_gu = &-\beta_0d_0\,r^{-\beta_0-2}\,f\,v +(n-1-
2\beta_0)r^{-\beta_0-1}\,f'\,v \cr &+
(n-1-2\beta_0)\,r^{-\beta_0-1}\,f\,\partial_rv + 2
r^{-\beta_0}\,f'\,\partial_rv + r^{-\beta_0}\,f''\,v \cr &+
r^{-\beta_0}\,f\,\partial_r^2v + r^{-\beta_0-2}\,f\,\Delta_\omega v + {\cal
L}(r^{-\beta_0}\,f\,v).}$$
We may thus write $\Delta_gu = u^q + Su$ in terms of the functions $f$ and $v$
as
$$\eqalign{&-\beta_0d_0\,v + (n-1-2\beta_0){rf'\over f}\,v +
(n-1-2\beta_0)\,r\,\partial_rv + 2{r^2f' \over f}\partial_rv \cr  &+ {r^2 f''
\over f}\,v + r^2\partial_r^2v + \Delta_\omega v = f^{q-1}\,v^q + {\cal
O}(r),}$$
since the derivative estimates and the form of the perturbation
operator ${\cal L}$ imply
$${r^{\beta_0q} \over f}\,{\cal L}(r^{-\beta_0}\,f\,v) = {\cal O}(r).$$
Using the change of variables $r=e^{-t}$ and applying the one-dimensional
version of the asymptotic maximum principle along with estimates on $f'$, $f''$,
and $v$ then allows us to conclude that $\overline{v}(\omega)$ weakly satisfies
$\Delta_\omega v \geq \beta_0d_0\,v$ on ${\bf S}^{n-1}\backslash L_p\Gamma.$
(Repeat the arguments in the proof of Proposition 2 in Section 5, and consult
[FM] for the details of the change of variables.)

\newpar
To conclude $f(r) \leq C\,r^{\epsilon}$, we compare $\overline{v}$ to a solution
of the linear equation $\Delta_\omega v = \beta_0d_0 v$ on ${\bf
S}^{n-1}\backslash L_p\Gamma$ that has a nonremovable singularity at
$L_p\Gamma$.  We prove the existence of such a solution by using the Poisson
transform developed in [Fn1].   Let $\{S_\alpha\}$ be a stratification of
$L_p\Gamma$ (exists since $T_p\Gamma$ is a proper tangent cone), and let
$\varphi_\alpha \in C(S_\alpha)\cap L^1(S_\alpha)$ be positive on $S_\alpha$.
(We use the induced metric on $S_\alpha$ to define a measure $\mu_\alpha$ on
$S_\alpha$.)  We now construct a solution on $S^{n-1}\backslash L_p\Gamma$  (or
a measure on $L_p\Gamma$ - these are equivalent cf.\ [Mz]) by using
the linearity of the equation to add solutions
$$v_\alpha(x) = \int\limits_{\overline {S_\alpha}} {\cal
G}(x,y)\,\varphi_\alpha(y)\,d\mu_\alpha(y)\leqno(6.9)$$
to $\Delta_\omega v = \beta_0d_0\,v$ on ${\bf S}^{n-1}\backslash S_\alpha$ where
${\cal G}(x,y)$ is the Green's function for $-\Delta_\omega+ \beta_0d_0$. (The
positivity of $\varphi_\alpha$ ensures that $v_\alpha$ is positive.) A simple
estimation of the integral in (6.9) (cf.\ [Fn1]) shows that as $x \to {\rm
int}(S_\alpha)$ that
$$\sigma(x)^{n-k-3}\,v_\alpha(x) \geq \epsilon \quad\hbox{if}\quad k < n-3$$
or
$$\big(1+\big|\log(\sigma(x)\big|\big)^{-1}\,v_\alpha(x) \geq \epsilon
\quad\hbox{if}\quad k = n-3$$
where $k$ is the dimension of $S_\alpha$.  Notice that the
assumption $\dim(T_p\Gamma) = \dim(L_p\Gamma) + 1 < d_0$ implies $k < n-3$.

\newpar
We now claim there exists a positive solution $v$ of $(\ast)$ such that
$\overline v \geq v$ on a path tending towards $L_p\Gamma$.  By the linearity of
the equation, we can clearly choose a positive solution $v$ of $(\ast)$
such that
for a given point $w_0 \in {\bf S}^{n-1}\backslash L_p\Gamma$ we have
${\overline v}(\omega_0) > v(\omega_0)$.  The difference function $w =
\overline v - v$ then satisfies $w(\omega_0) > 0$ and in the set $\Omega =
\{\omega \in {\bf S}^{n-1}\backslash L_p\Gamma: w > 0\}$ we have $\Delta_\omega
w \geq 0$, so  the maximum principle then implies that $\sup_{\partial \Omega} w
> 0$.  We can therefore conclude that there exists a path on ${\bf
S}^{n-1}\backslash L_p\Gamma$ tending to $L_p\Gamma$ with $\overline v\geq
v$, and we may also assume that $\overline v \geq \delta\,\sigma^{-(n-k-3)}$ for
some positive constant $\delta$ where $k= \dim(L_p\Gamma)$.  Hence, $u(r,\omega)
\geq \delta\,r^{-\beta_0}\,f(r)\,\sigma^{-(n-k-3)}$ close to $p$.  Choosing
$\sigma$ proportional to $r$, and using the fact that $\rho \approx \sigma\,r$
(see the proof of Proposition 3 in Section 5), we find that
$$u \geq C\,r^{-\beta_0-(n-k-3)}\,f(r)$$
on a path in $M\backslash\Gamma$ approaching $p$.  But, from the upper bound in
Lemma 1, we know that $u \leq C\rho^{-\beta_0}$ on this path which implies that
$f(r)$ must satisfy
$$f(r)\,r^{-\beta_0-(n-k-3)} \leq C\,r^{-2\beta_0},$$
or
$$f(r) \leq C\,r^{-\beta_0 + n - k - 3} \leq C\,r^{d_0 - k - 1}.$$
The assumption that $\dim(T_p\Gamma) < d_0$ then implies
that $k < d_0 -1$ and thus $f(r) \leq C\,r^\epsilon$.  From which,
a simple calculation shows $\hat g = u^{q-1}\,g$ is not complete. \qed

\newpar\noindent
{\bf Remark:} The proof breaks down exactly when $k =d_0-1$.  In this case, we
can conclude only that $f(r) \leq C$. A possibly cause why our proof breaks
down in this case is that we do not strongly use the ``negative case'' in this
last part of the proof; in the positive case (cf.\ [Pa]) when $q={n+2
\over n-2}$ there exist solutions when $\dim(\Gamma) = d_0$ where the asymptotic
growth is on the order of
$$u \approx r^{-(n-2)/2}\,\log(1/r)^{-(n-2)/4}.$$
And moreover, this type of growth implies the conformal metric $\hat g =
u^{4/(n-2)}\,g$ is a complete metric.

\vfill\eject

\centerline {\bf Bibliography}

\ref {Au} {T.\ Aubin, {\it Nonlinear Analysis on Manifolds},
Springer-Verlag,  (1982).}

\ref {Av} {P.\ Aviles,{\it A study of the singularities of solutions of a
class of nonlinear elliptic partial differential equations},  Comm.\
PDE, {\bf 7} (1982) 609-643.}

\ref {AM1} {P.\ Aviles and R.\ McOwen, {\it Conformal deformation to
constant negative scalar curvature in complete Riemannian manifolds}, J.\
Diff.\  Geom.,  {\bf 27} (1988) 225-238.}

\ref {AM2} {P.\ Aviles and R.\ McOwen, {\it Complete conformal metrics with
negative scalar curvature in compact Riemannian manifolds}, Duke Math.\ J.,
{\bf 56} (1988) 395-397.}

\ref {BP} {P.\ Baras and M.\ Pierre, {\it Singularit\'es \'eliminables
pour des equations semi-lin\'eaires},  Ann.\ Inst.\ Fourier Grenoble, {\bf
34I} (1984) 185-206.}

\ref {Ca} {E.\ Calabi, {\it An extension of E. Hopf's maximum principle
with  an application to Riemannian geometry}, Duke Math.\ J., {\bf 25}
(1958) 45-56.}

\ref {CY} {S.Y.\ Cheng and S.T.\ Yau, {\it Differential equations on
Riemannian manifolds and their geometric applications}, Comm.\ Pure.\
Applied.\ Math {\bf XXVIII} (1975) 333-355.}

\ref {De} {P.\ Delano\"e, ``Generalized stereographic projections with
prescribed scalar curvature", {\it Contemporary Mathematics: Geometry,
Physics, and Nonlinear PDE}, ed.\ V.\ Oliker and A.\
Treibergs, AMS (1990).}

\ref {Fe} {H.\ Federer, {\it Geometric Measure Theory}, Springer-Verlag
(1969)}

\ref {Fn1} {D.\ Finn, {\it Positive solutions of $\Delta_gu = u^q + Su$
singular at submanifolds with boundary}, Indiana Univ.\ Math.\ J.\ {\bf 43}
(1994) 1359-1397.}

\ref {Fn2} {D.\ Finn, {\it Positive solutions to nonlinear elliptic equation
with prescribed singularities}, PhD Thesis, Northeastern University, 1995}

\ref {FM} {D.\ Finn and R.\ McOwen, {\it Singularities and asymptotics for
the equation $\Delta_gu - u^q =Su$}, Indiana Univ. Math. J., {\bf 42} (1993)
1487-1523.}

\ref {GT} {D.\ Gilbarg and N.\ Trudinger, {\it Elliptic Partial
Differential Equations of Second Order}, Springer-Verlag (1983).}

\ref {GM} {M.\ Goresky and R.\ MacPherson, {\it Stratified Morse Theory},
Springer-Verlag (1988).}

\ref {Kz} {J.\ Kazdan {\it Prescribing the Curvature of a Riemannian
Manifold}, AMS (1985).}

\ref {LP} {J.\ Lee and T.\ Parker, {\it The Yamabe problem}, Bull.\ Amer.\
Math.\ Soc., {\bf 17} (1987) 37-91.}

\ref {LTY} {P.\ Li, L.\ Tam, and D.\ Yang, {\it On the elliptic equation
$\Delta u + ku - K u^p = 0$ on complete Riemannian manifolds and their
geometric applications I}, preprint.}

\ref {LN} {C.\ Loewner and L.\ Nirenberg, ``Partial differential equations
invariant under conformal or projective transformations", {\it Contributions
to Analysis}, Academic Press, 245-275, (1975).}

\ref {Mz} {R.\ Mazzeo, {\it Regularity for the singular Yamabe problem},
Indiana Univ.\ Math.\ J. {\bf 40} (1991) 1227-1299.}

\ref {MP} {R.\ Mazzeo and F.\ Pacard, {\it A construction of
singular solutions for a semilinear equation using asymptotic analysis},
(to appear in J. Diff. Geom).}

\ref {MPU} {R.\ Mazzeo, D.\ Pollack, K.\ Uhlenbeck, {\it Moduli spaces of
singular Yamabe metrics}, (to appear in Jour.\ of AMS)}

\ref {Mc1} {R.\ McOwen, ``Singularities and the conformal scalar curvature
equation'', {\it Geometric Analysis and Nonlinear PDE}, ed I.\ Bakelma, Marcel
Dekker (1992).}

\ref {Mc2} {R.\ McOwen, {\it Partial Differential Equations: Methods and
Applications}, Prentice Hall, 1996.}

\ref {Mi} {J.\ Milnor, {\it Singular Points of Complex Hypersurfaces},
Princeton University Press (1968).}

\ref {P} {F.\ Pacard, {\it The Yamabe problem on subdomains of even
dimensional spheres}, (to appear in J.\ Diff.\ Geom.).}

\ref {R} {H.L.\ Royden, {\it Real Analysis}, Macmillan, (1988)}

\ref {Sc1} {R.\ Schoen, {\it Conformal deformation of a Riemannian metric to
constant scalar curvature}, J.\ Diff.\ Geom, {\bf 20}, (1984) 479-495}

\ref {Sc2} {R.\ Schoen, {\it The existence of weak solutions with
prescribed singular behavior for a conformally invariant scalar equation},
Comm.\ Pure.\ Applied.\ Math.,{\bf XLI} (1988) 317-392.}

\ref {SY1} {R.\ Schoen and S.T.\ Yau, {\it On the structure of manifolds
with positive curvature},  Manu.\ Math., {\bf 28}, (1979) 159-183.}

\ref {SY2} {R.\ Schoen and S.T.\ Yau, {\it Conformally flat manifolds,
Kleinian groups and scalar curvature},  Invent.\ Math.\ {\bf 92} (1988)
47-71.}

\ref {Th} {R.\ Thom, {\it Ensembles et Morphisms Stratifi\'es}, Bull.\
Amer.\ Math.\ Soc.\ {\bf 75} (1969) 240-284}

\ref {Vn} {L.\ Veron, {\it Singularit\'es \'eliminables d'\'equations
elliptiques non lin\'eaires}, J.\ Diff.\ Eq., {\bf 41} (1981) 87-95.}

\ref {Wh1} {H.\ Whitney, ``Local properties of analytic varieties'', {\it
Differential and Combinatorial Topology}, ed.\ S.\ Cairns, Princeton Univ.\
Press, (1965) 205-244.}

\ref {Wh2} {H.\ Whitney, {\it Tangents to an analytic variety}, Ann.\ Math
{\bf 81}, (1965) 496-549.}

\bye